\documentclass[12pt]{amsart}
\usepackage[utf8]{inputenc}

\usepackage[margin=1in]{geometry}
\usepackage{amsmath,amssymb,multicol,cancel,xcolor,graphicx,mathrsfs,amsthm,tikz}
\usepackage[shortlabels]{enumitem}
\usepackage{comment}
\usepackage{thm-restate}
\usepackage{hyperref}
\usepackage{chemfig}
\usepackage[utf8]{inputenc} 
\usepackage[backend=biber,style=alphabetic,maxnames=99,maxalphanames=6,giveninits=true]{biblatex}
\hypersetup{
	colorlinks=true, 
	linktoc=all,     
	linkcolor=black,  
}
\usepackage{subcaption}

\usepackage{tcolorbox}

\newcommand{\vertiii}[1]{{\left\vert\kern-0.25ex\left\vert\kern-0.25ex\left\vert #1 
		\right\vert\kern-0.25ex\right\vert\kern-0.25ex\right\vert}}
\newcommand{\R}{\mathbb{R}}

\newcommand{\Z}{\mathbb{Z}}

\newcommand{\hsp}{\hspace{.5cm}}

\newcommand{\inv}[1]{#1^{-1}}

\newcommand{\mc}[1]{\mathcal{#1}}

\newcommand{\N}{\mathbb{N}}

\newcommand{\bd}{\partial}

\newcommand{\Q}{\mathbb{Q}}

\newcommand{\ind}{\operatorname{ind}}

\numberwithin{equation}{section}
\newtheorem{thm}{Theorem}[section]
\newtheorem{lemma}[thm]{Lemma}
\newtheorem{prop}[thm]{Proposition}
\newtheorem{cor}[thm]{Corollary}

\newtheorem{conj}[thm]{Conjecture}

\newtheoremstyle{uprightnormal}
{3pt}{3pt}
{\normalfont}    
{}{\bfseries}
{.}{.5em}{}

\theoremstyle{uprightnormal}
\newtheorem{defn}[thm]{Definition}
\newtheorem{rmk}[thm]{Remark}
\newtheorem{nota}[thm]{Notation}
\newtheorem{examp}[thm]{Example}

\title{Filtered Instantons and the Concordance of Satellites}
\author{Ivan So}
\address{Department of Mathematics, Michigan State University, East Lansing, MI 48824}
\email{soivan@msu.edu}

\bibliography{reference.bib}
\begin{document}
	\small
	\maketitle
	\begin{abstract}
		In \cite{NST23}, Nozaki-Sato-Taniguchi defined a family of invariants $ r_s $ for integer homology spheres with filtered instanton homology \cite{FS92}. Coupling these with techniques in classical knot theory, we produce some results in the knot concordance group, including criteria for a family of satellite knots to be linearly independent and the independence of what we call the Whitehead $ n $-ble $ P_n $.
	\end{abstract}
\section{Introduction}
Two knots $K_0,K_1\subset S^3$ are said to be \textit{concordant} if there exists an embedding $f:S^1\times I\hookrightarrow S^3\times I$ such that the restriction of the embedding to the two boundary components are $ K_0 $ and $ K_1 $. Concordance is clearly an equivalence relation. Indeed, by taking the monoid of connected sum of knots quotienting by the concordance, an abelian group is defined by taking the identity to be the concordance class of the unknot $ U $. The abelian group defined under this construction is known as the \textit{concordance group}, denoted by $\mathcal{C}$. In the literature, the concordance group is studied in various categories, including the topological and smooth categories. In this paper, we focus on the \textit{smooth} concordance group.

Let $ P\subset S^1\times D^2 =V$ be a knot. The \textit{satellite knot} $ P(K)\subset S^3 $ for a knot $ K\subset S^3 $ is defined as the image of $ P $ under the embedding of $ V $ into a tubular neighborhood of $ K $ with 0-framing. In this context, $ P $ is known as the \textit{pattern} or \textit{satellite operator} and $ K $ is known as the \textit{companion knot}. The behavior of $ \mc{C} $ under the satellite operator is a frequently considered question in knot theory. 

Concordance is closely related to the study of homology cobordism group of $ \Z/2 $-homology 3-spheres $ \Theta^3_{\Z/2} $ as branched double cover of knot defines a well-defined homomorphism $ \Sigma_2:\mc{C}\to\Theta^3_{\Z/2} $. 

Instanton Floer homology is particularly strong in distinguishing integer homology 3-spheres constructed from $ 1/n $-surgery over the same knot \cite{Fur90}, \cite{FS90}. In this paper, we exploit this property to study the smooth concordance group under certain satellite operation using a numerical invariant $r_s$ \cite{NST23} arising from the filtered instanton Floer homology \cite{FS92}. 
\subsection{Independence in the concordance group}
The most general form of our result is as follows:
\begin{thm}\label{main}
	Let $Q_n\subset S^1\times D^2$ be a family of patterns such that
	\begin{enumerate}[(1)]
		\item 
		$ Q_n $ is rational unknotting number 1 within $ S^1\times D^2 $;
		\item 
		$\Sigma_2(Q_n(U))\cong S^3_{1/n}(J)$ and;
		\item 
		$h(S^3_1(K\# J\# K))<0$,
	\end{enumerate}   
	where $h:\Theta^3_\Z\to\Z$ is the Fr\o yshov homomorphism from \cite{Fro02}, then $\{Q_n(K)\}_{n\in\mathbb{N}}$ are linearly independent in $\mc{C}$. 
\end{thm}
Condition (3) in the above theorem is difficult to apply as the map $ h(S^3_1(-)):\mc{C}\to\Z $ is not a homomorphism and and is generally hard to compute. To this end, \cite{DISST22} defined a \textit{slice torus invariant} $ \widetilde{s} $ whose positivity suffices to conclude that  $ h(S^3_1(K))<0 $. Hence, Theorem \ref{main} has the following hands-on corollary
\begin{cor}\label{maintorus}
	Let $Q_n\subset S^1\times D^2$ be a family of patterns such that
	\begin{enumerate}[(1)]
		\item 
		$ Q_n $ is rational unknotting number 1 and;
		\item 
		$\Sigma_2(Q_n(U))\cong S^3_{1/n}(J)$ for a knot $J$ such that $2\widetilde{s}(K)+\widetilde{s}(J)>0$,
	\end{enumerate} 
	then $\{Q_n(K)\}_{n\in\mathbb{N}}$  are linearly independent in $\mc{C}$.
	
	In particular, if $K,J$ are alternating knots and $2\sigma(K)+\sigma(J)<0$, where $\sigma$ is the Seifert signature of a knot, then $\{Q_n(K)\}_{n\in\N}$ are linearly independent in $\mathcal{C}$.
\end{cor}

\begin{figure*}[h]
	\centering

	\begin{subfigure}[t]{0.5\textwidth}
		\centering
		\includegraphics[height=1.2in]{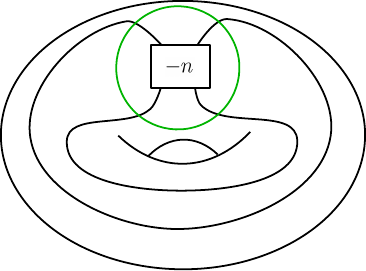}
		\caption{The pattern $ P_n $}
	\end{subfigure}
	~
	\begin{subfigure}[t]{0.5\textwidth}
		\centering
		\includegraphics[height=1.2in]{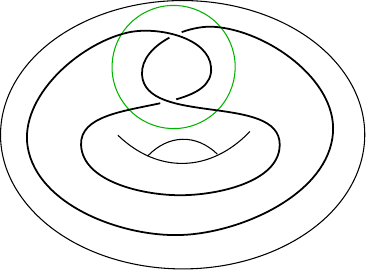}
		\caption{$ P_2=Wh_+ $, the untwisted positive-clasped Whitehead pattern}
	\end{subfigure}%
	
	\caption{The pattern $ P_n $ considered in this paper. Here $ n $ denote the number of positive half-twist.}
	
	\label{whiteheadlike}
\end{figure*}

If we restrict ourselves to the Whitehead $n$-ble pattern $ P_n $ (see Figure \ref{whiteheadlike}), then the algorithm introduced in Section 2 will imply $ J=U $. As a corollary

\begin{cor}\label{whiteheadnble}
	Let $ P_n $ be the Whitehead $ n $-ble pattern. Then $ \{P_n(T_{p,q})\}_{n\in\mathbb{N}}$ are linearly independent in $ \mc{C} $ for any fixed coprime $ p,q>1 $. 
\end{cor}

In the proof of Corollary \ref{whiteheadnble}, we will show the corresponding value of $ r_s $ for the branched covers is finite. Since $ r_s(S^3)=\infty $, this implies a slightly stronger result than that of Pinz\'on-Caicedo in \cite{Pin17} (which applies only for even $ n $)

\begin{cor}\label{whitenotslice}
	$ \{P_n(T_{p,q})\}_{n\in\mathbb{N}}$ are not slice for any fixed coprime $ p,q>1 $. 
\end{cor}

The computation of the Fr\o yshov $ h $-invariant is never easy. However, due to this difficulty, we expect the  condition in Theorem \ref{main} will provide a lot of possibilities for linear independence-type results provided we can find ways to compute $ h $ for certain 1-surgeries.


\subsection{Linearly independence under fixed Whitehead $ n $-ble}

The concordance group $ \mc{C} $ is an abelian group. In this sense, the pattern $ P\subset S^1\times D^2 $ defines a satellite operator $ P:\mc{C}\to\mc{C} $. One of the questions in knot theory is to ask whether satellite operators generate a ``large'' images subgroup. Following \cite{HPC21}, we define
\begin{defn}
	A pattern $P\subset S^1\times D^2$ is said to be of \textit{infinite rank} if there exists a family of knots $\{K_n\}_{n\in\N}$ such that $\{P(K_n)\}_{n\in\N}$ gives a $\Z^\infty$-subgroup of $\mc{C}$.
\end{defn}

An example of result of this type  is as follows

\begin{thm}[\cite{HK12}, Theorem 1]
	$ \{P_2(T_{2,2^n-1})\} $ for $ n=2,3,\cdots $ are linearly independent in $ \mc{C} $.
\end{thm}

Hedden-Pinzón-Caicedo \cite{HPC21} conjectured that all non-constant satellite operators have infinite rank. Beyond asking whether a particular pattern has infinite rank, a natural question is: \textit{which} companion knots yield an infinite-rank image? We have the following result:

\begin{thm}\label{almostrankexpand}
	Let $M=\{(m_k,k)\in\N\times\N|m_k\leq m_{k+1}\}$. Fix $ n\in\N $ and coprime $ p,q>1 $,
	then $ \{P_n(m_kT_{p,q+kp})\}_{(m_k,k)\in M} $ are linearly independent in $ \mc{C} $. 
\end{thm}

\begin{rmk}
	Theorem 13 of \cite{Pin17}  is similar to Theorem \ref{almostrankexpand} and Corollary \ref{whiteheadnble}, but in a slightly different direction. Their result requires the set $ \{(n_i,(p_i,q_i))\}_{i\in\N} $ satisfies $$ p_iq_i(2n_ip_iq_i)<p_{i+1}q_{i+1}(n_{i+1}p_{i+1}q_{i+1}-1). $$ While our results do not allow simultaneous variation of $ (p,q) $ and $ n $, they do cover cases that were inconclusive in prior work. For example,  for $ (p,q)=(2,3) $ and $ n=2,4 $, \cite{Pin17} is inconclusive on the independence whereas independence follows from our computation. Corollary 1.12 of \cite{DISST22}  is also similar in nature, but our result cover a larger class of companion knots.
	
	In the case of $ P_{2n-1} $,  Hedden-Pinzón-Caicedo (\cite{HPC21}, Proposition 8) showed all patterns with nonzero winding number are necessarily infinite rank. Therefore, the above theorem does not provide a new infinite-rank result. Their proof, however, was non-constructive. The above theorem also gives an explicit family of knots whose images form an infinite-rank subgroup for the family of patterns $ \{P_{2n-1}\}_{n\in\N} $.
\end{rmk}

\subsection{Comparison to prior technologies} Prior to our work, Hedden-Kirk \cite{HK12} and Pinzón-Caicedo \cite{Pin17} proved related results using $ SO(3) $ instanton theory. Their linear independence results in the concordance group relied on the computation of the minimal value Chern-Simons invariant. In contrast, the invariant $ r_s $ is an invariant from $ SU(2) $ instanton theory extracted from certain element over instanton Floer cohomology group, which contains more information than just the minimal value Chern-Simons invariant. we believe this is the primary reason why a somewhat stronger independence result can be extracted.

Compared to Heegaard Floer homology, instanton Floer homology is certainly more difficult to compute and thus remains inconclusive for a large class of satellite operator. But to the best of the author's knowledge, independence results for $ 1/n$-surgeries on a fixed knot remain a unique feature for instanton Floer homology (see for instance \cite{Fur90}, \cite{DLM24}).

\noindent \textbf{Organization}: In Section 2, we describe an algorithm for obtaining a surgery description of the branched cover of satellite knots. In Section 3, we provide an overview of the instanton filtration invariant $ r_s $, including its definition and key properties. The proofs of the main theorems are contained in Section 4. In Section 5, we present an application involving a nontrivial Montesinos pattern and a non-torus companion knot. \\

\noindent \textbf{Acknowledgement}: The author would like to thank his advisor, Matthew Hedden, for helpful discussions and comments on an earlier draft of this paper; Masaki Taniguchi for pointing out his earlier work which is closely related to our results; and Jennifer Hom for some interesting discussions.

\section{Surgery Description of Branched Double Cover}
In this section, we recall the algorithm presented in Section 2 of \cite{DHMS24} which identifies branched double cover of certain satellite knots in terms of surgery.

Such an algorithm is helpful for computing Floer-type invariants, as these invariants tend to be difficult to calculate directly from the branched cover description. However, with a surgery description, one can use tools such as the surgery exact triangle or surgery formulas to compute Floer homologies and hence the associated numerical invariants.

For example, in the recent notable disprove of the sliceness of $ (2,1) $-cable of figure eight knot $ 4_1 $ (denoted by $ (4_1)_{2,1} $), \cite{DKMPS24} used the fact that 
$$ \Sigma_2((4_1)_{2,1})\cong S^3_1(4_1\# 4_1).$$ 
As a result, the calculation in Heegaard-Floer theory reduces to manageable problem: the knot Floer complex of $ 4_1 $ is known, there is a connected sum formula for knot Floer complex and the Heegaard-Floer complex of $ S^3_1(4_1\# 4_1) $ can be computed via the integer surgery formula.

To describe the algorithm, we need the notion of \textit{rational tangles}.	
\subsection{Rational tangles}
\textbf{add more details about rational tangles.} We begin with the definition of rational tangles. In the following, we are always dealing with tangles in $ B^3 $ with four ends on $ \bd B^3 $.
\begin{defn}
	A tangle is said to be \textit{rational} if it consists of a pair of boundary-parallel arcs.
\end{defn}

\begin{figure*}[h]
	\begin{subfigure}[t]{0.25\textwidth}
		\centering
		\includegraphics[scale=0.8]{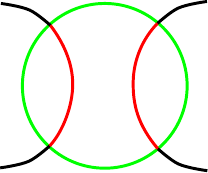}
		\caption{ $ T_{\infty} $}
	\end{subfigure}
	~
	\begin{subfigure}[t]{0.25\textwidth}
		\centering
		\includegraphics[scale=0.8]{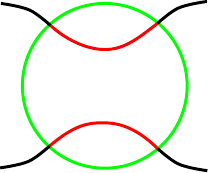}
		\caption{$ T_0 $}
	\end{subfigure}
	~
	\begin{subfigure}[t]{0.25\textwidth}
		\centering
		\includegraphics[scale=0.8]{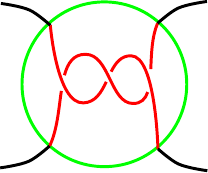}
		\caption{$ T_{3} $}
	\end{subfigure}
	~
	\begin{subfigure}[t]{0.25\textwidth}
		\centering
		\includegraphics[scale=0.8]{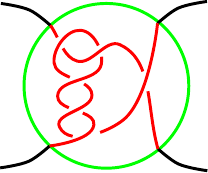}
		\caption{$ T_{[-1,3,-2]}=T_{-3/5} $}
	\end{subfigure}
	\caption{Example of rational tangles.}
	
	\label{rationaltangles}
\end{figure*}

\begin{examp}
	More explicitly, our tangle of interest will have four ends, which we denote as NW, NE, SW, and SE, corresponding to the four cardinal directions on a compass (these are the ends of the red arcs in Figure \ref{rationaltangles}). Let $h$ represent a positive half-twist with respect to the horizontal (E–W) axis, while $v$ represents a positive half-twist with respect to the vertical (N–S) axis.
	
	To build a rational tangle, start with a circle with four marked points (NW, NE, SW, SE). First, attach a tangle with $a_1$ positive half-twists, corresponding to the tangle word $h^{a_1}$, to NE and SE. If this suffices to describe the tangle, connect the two loose ends to NW and SW. Otherwise, treat these loose ends as the new NE and SE and attach vertical half-twists to SW and SE, corresponding to the tangle word $v^{a_2}$. By iterating this process, the resulting tangle word is of the form
	$$ v^{a_{2m}}h^{a_{2m-1}}\cdots v^{a_2}h^{a_1}. $$ 
	Figure \ref{rationaltangles} (3) and (4) provide examples of rational tangles enumerated by the above process.
	
	Let 
	$$[a_1,\cdots,a_n]:=a_1+\frac{1}{a_2+\frac{1}{\cdots+\frac{1}{a_n}}}.$$
	With $ T_0 $ and $ T_\infty $ defined as in (1) and (2) of Figure \ref{rationaltangles}, we label tangles by $T_{[a_1,\cdots,a_n]}$. Indeed, \cite{Con70} showed that the set of all rational tangles is in one-to-one correspondence with $ \Q\cup\{\infty\} $.
	
	Note that this enumeration is relative: one can choose different tangles as $ T_0 $ and $ T_\infty $, resulting in different enumerations than those presented in Figure \ref{rationaltangles}. This relativity is important for the computation of surgery coefficients in the Montesinos trick, which is discussed in later sections. Meanwhile, we call the enumeration with $ T_\infty $ and $ T_0 $ defined in Figure \ref{rationaltangles} the \textit{standard enumeration}.
	
\end{examp}

Rational tangles are mentioned for the following definition.
\begin{defn}
	Let $ P\subseteq S^1\times D^2 $ be a pattern. We say $ P $ has \textit{rational unknotting number one} if there exists a rational tangle $ T $ in $ P $ such that replacing $ T $ with another rational tangle $ T' $ gives an unknot in $ S^1\times D^2 $. 
\end{defn}

\subsection{The Montesinos trick for patterns of rational unknotting number one}

Given a pattern $ P\subset S^1\times D^2 $ of rational unknotting number one, let $ P' $ be the unknot in $ S^1\times D^2 $ obtained from rationally unknotting  $ P $. The Montesinos trick \cite{Mon75} implies $ \Sigma_2(P(U)) $ can be obtained from $ \Sigma_2(P'(U))\cong S^3 $ by surgery over certain \textit{strongly invertible} knot $ J $ (a knot with an involution $ \tau $ which fixes two points on the knot). We call such $ J $ the \textit{Montesinos friend} of $ P $.

In other words, $ \Sigma_2(P(U))\cong S^3_{p/q}(J) $ and the branched involution is identified with the involution on $ S^3_{p/q}(J) $ induced by the involution on $ J $.

In the following, we will assume the gluing of solid torus to knot complement to be the one that identifies the meridian of solid torus with the Seifert framing of $ K $ unless specified.

Given a pattern $ P\subset S^1\times D^2 $ and let $ \mu $ be its meridian, with a companion knot $ K $, $ P(K) $ can be taken as the image of $ P $ in the gluing
$$S^3\cong (S^3-N(\mu))\cup_{\bd N(\mu)}(S^3-N(K)).$$
Using the picture from Montesinos trick and let $ \widetilde{\mu} $ be a meridian of $ J $, we have
$$\Sigma_2(P(K))\cong (S^3_{p/q}(J)-N(\widetilde{\mu})-N(\tau\widetilde{\mu}))\cup_{\bd N(\widetilde{\mu})} (S^3-N(K))\cup_{\bd N(\tau\widetilde{\mu})}(S^3-N(K)).$$
In the surgery picture, this means
$$\Sigma_2(P(K))\cong S^3_{p/q}(J_{\widetilde{\mu},K}),$$
where $ J_{\widetilde{\mu},K} $ is the infection of $ J $ by $ K $ at the meridians $ \widetilde{\mu} $ and $ \tau\widetilde{\mu} $ (in case infecting over the meridians, $ J_{\widetilde{\mu},K} $ is just $ K\#J\# K$).

\begin{figure*}[h]
	\begin{subfigure}[t]{0.3\textwidth}
		\centering
		\includegraphics[height=1.2in]{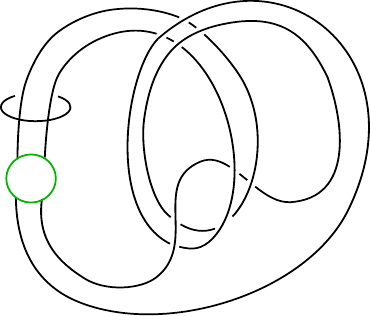}
		\caption{}
	\end{subfigure}
	~
	\begin{subfigure}[t]{0.3\textwidth}
		\centering
		\includegraphics[height=1.2in]{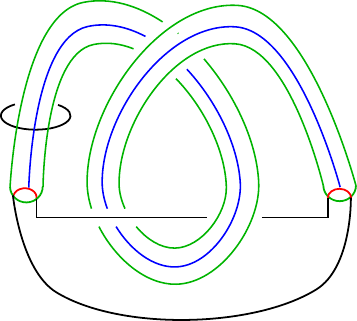}
		\caption{}
	\end{subfigure}
	~
	\begin{subfigure}[t]{0.3\textwidth}
		\centering
		\includegraphics[height=1.2in]{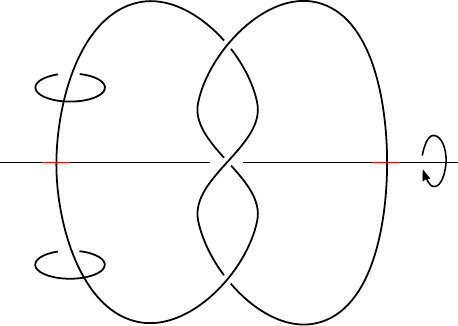}
		\caption{}
	\end{subfigure}
	\caption{Realization of $ J $.}
	
	\label{realization}
\end{figure*}

\begin{figure*}[h]
	\begin{subfigure}[t]{0.3\textwidth}
		\centering
		\includegraphics[height=1.2in]{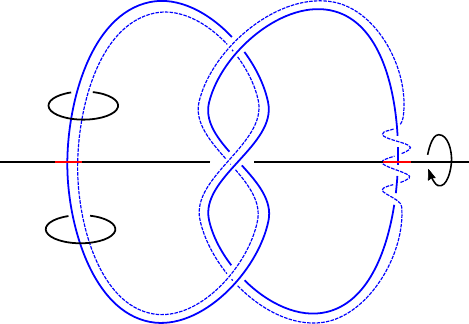}
		\caption{}
	\end{subfigure}
	~
	\begin{subfigure}[t]{0.3\textwidth}
		\centering
		\includegraphics[height=1.2in]{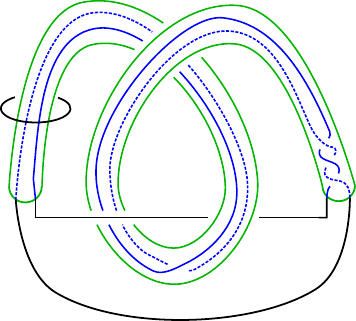}
		\caption{}
	\end{subfigure}
	~
	\begin{subfigure}[t]{0.3\textwidth}
		\centering
		\includegraphics[height=1.2in]{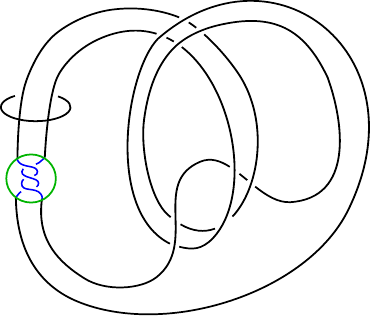}
		\caption{}
	\end{subfigure}
	\caption{Computation of the surgery coefficient.}
	
	\label{surgerycoefficient}
\end{figure*}

From the prior discussions, finding the surgery description of $ \Sigma_2(P(K)) $ for a given pattern $ P $ with rational unknotting number one is equivalent to answering the following questions:
\begin{enumerate}
	\item 
	What is the Montesinos friend $ J $ of $ P $?
	\item 
	What is the surgery coefficient $ p/q $?
\end{enumerate}
Section 2 of \cite{DHMS24} described an algorithm that answers both questions. We illustrate their algorithm with their illuminating example presented in Figure \ref{realization} and \ref{surgerycoefficient} (borrowed from Figure 4 of \cite{DHMS24}).

To find $ J $, first remove a rational tangle to be unknotted and replace it by $ T_\infty $ as illustrated in (1), (2) of Figure \ref{realization}. Isotope the $ B^3 $ containing replaced tangle such that remaining part of the knot (in black) and the newly replaced tangle (in red) form an unknot. The core of the cylinder (in blue) is the quotient of $ J $ under the branched involution. The knot $ J $ will be the lift of the core of the solid cylinder, given in (3) of Figure \ref{realization}.

To find the surgery framing $ p/q $, we need the $ \tau $-equivariant Seifert framing of $ J $. The $ \tau $-equivariant Seifert framing of the example considered in Figure 3 is presented in (1) of Figure \ref{surgerycoefficient}. (2) of Figure \ref{surgerycoefficient} shows the quotient of the knot together with the framing under the branched involution over $ U\subset S^3 $. Isotope the green circle which represents the $ B^3 $ containing the tangle back to the standard $ B^3 $, the rational number corresponding to the rational tangle (in the standard enumeration) is our \textit{reference framing} for a pattern with Montesinos friend $ J $.

Let $ T_{p'/q'} $ be the rational tangle (in the standard enumeration) replaced under rational unknotting and suppose $ r/s $ is the reference framing for $ J $. The desired surgery coefficient $ p/q $ in the Montesinos trick is given by $p/q=p'/q'-r/s$. (e.g. if the tangle to be replaced in the green circle of (1) of Figure \ref{realization} is the tangle $ p'/q' $, then the surgery coefficient $ p/q=p'/q'-3 $).

\subsection{Surgery description of $ \Sigma_2(P_n(K)) $}
\begin{figure*}[h]
	\begin{subfigure}[t]{0.3\textwidth}
		\centering
		\includegraphics[height=1.2in]{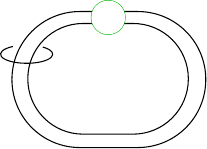}
		\caption{}
	\end{subfigure}
	~
	\begin{subfigure}[t]{0.3\textwidth}
		\centering
		\includegraphics[height=1.2in]{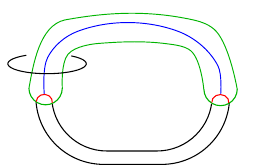}
		\caption{}
	\end{subfigure}
	~
	\begin{subfigure}[t]{0.3\textwidth}
		\centering
		\includegraphics[height=1.2in]{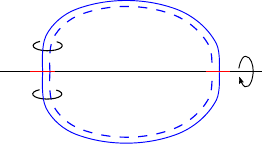}
		\caption{}
	\end{subfigure}
	\caption{Determination of $ J $ for the pattern $ P_n $.}
	
	\label{Pn}
\end{figure*}

Let $ P_n\subset S^1\times D^2 $ be the knots presented in Figure 1. By removing the rational tangle encircled in green, this gives us the diagram in (1) of Figure \ref{Pn}. In (2) of Figure \ref{Pn}, an isotopy after the replacement with $ T_\infty $. The blue arc traces the quotient of $ J $ under the branched involution of $ \Sigma_2(U) $ is shown. (3) in Figure \ref{Pn} then shows $ J=U $  with the $ \tau $-invariant framing.

According to the algorithm in the previous subsection, the surgery coefficient  $ p/q $ depends only on the rational tangle encircled in Figure 1, which coresponds to $ \frac{1}{n}\in\Q\cup\{\infty\} $. As in the case $ J=U $ , $ J_{\widetilde{\mu},K}=K\# K $. This implies:
\begin{lemma}\label{Montesinos}
	$ \Sigma_2(P_n(K))\cong S^3_{1/n}(U_{\widetilde{\mu},K})\cong S^3_{1/n}(K\# K)$.\qed
\end{lemma}

\section{Filtered instanton Floer homology and the instantons filtration invariant}
\subsection{Filtered instantons and the definition of $r_s$}
In this subsection, we recall the numerical invariant for integer homology 3-sphere defined in \cite{NST23}.
We assume the readers are familiar with the basics of instanton Floer homology \cite{Flo88}, \cite{Don02}. The following notations will be used throughout this section:
\begin{nota}
	\;
	\begin{itemize}
		\item
		$ \mc{A}(Y) $, the set of $ SU(2) $-connections on the product $ SU(2) $-bundle.
		\item 
		${\mc{B}}(Y):=\mc{A}(Y)/\operatorname{Map}_0(Y,SU(2))$.
		\item 
		$ \widetilde{\mc{B}}(Y):=\mc{A}(Y)/\operatorname{Map}_0(Y,SU(2)) $, where $\operatorname{Map}_0(Y,SU(2))$ denote the set of smooth functions of mapping degree 0.
		\item 
		$ \mc{B}^*(Y) \subset\mc{B}(Y)$, the subspace consists only of irreducible connections.
		\item 
		$ cs_Y:\widetilde{\mc{B}}(Y)\to\R $, the Chern-Simons functional.
		\item 
		$ \widetilde{R}(Y)\subset \widetilde{\mc{B}}(Y)$, the set of gauge equivalence class of flat connections on $ Y $ or equivalently, the set of critical points of the Chern-Simons functional ($ \widetilde{R}^*(Y) $ means irreducible flat connections, respectively).
		\item 
		$ \Lambda_Y:=cs_Y(\widetilde{R}(Y)) $, $ \Lambda^*_Y:=cs_Y(\widetilde{R}^*(Y)) $.
		\item 
		$ \R_Y:=\R\setminus\Lambda_Y $.
	\end{itemize}
\end{nota}
We first recall the notion of filtered instanton Floer homology by Fintushel-Stern in \cite{FS92}.
\begin{defn}
	Fix $ s\in[-\infty,0] $. For a given $ r\in\R_Y $, metric $ g $ on $ Y $, a non-degenerate regular perturbation $ \pi\in\mc{P}(Y,r,s,g) $ and orientations on line bundle $L_X$ for a 4-manifold $ X $ with $ \bd X=-Y $, the \textit{chains of the filtered instanton Floer homologies} are defined by 
	$$CI^{[s,r]}_i(Y):=\begin{cases}
		\Z\{[a]\in\widetilde{R}^*(Y)_\pi|\ind(a)=i,\; s<cs_{Y,\pi}(a)<r\} & \text{if }s\in\R_{Y}\\
		\Z\{[a]\in\widetilde{R}^*(Y)_\pi|\ind(a)=i,\; s-\lambda_Y/2<cs_{Y,\pi}(a)<r\} & \text{if }s\in\Lambda_{Y}
	\end{cases},$$
	where $ \lambda_Y:=\min\{|a-b|\;|a\neq b,\;a,b\in\Lambda_Y\} $.\\
	The \textit{boundary map} $ \bd^{[s,r]}:C_i^{[s,r]}(Y)\to C_{i-1}^{[s,r]}(Y) $ is given by the restriction of $ \bd $ to $ CI^{[s,r]}_i(Y,\pi) $.
\end{defn}

\begin{defn}
	The \textit{filtered instanton Floer homology} is defined as
	$$I^{[r,s]}_*(Y):=\ker\bd^{[s,r]}/\operatorname{im}\bd^{[s,r]}.$$
\end{defn}

In \cite{Tan18}, Taniguchi defined an obstruction invariant $ [\theta_Y]\in I^1(Y,\pi) $ for oriented homology spheres in the classical instanton Floer homology. The obstruction class can be generalized to the context of filtered instantons:
\begin{thm}[\cite{Tan18}, generalized to filtered instantons]
	The homomorphism
	\begin{align*}
		\theta^{[s,r]}_Y:CI^{[s,r]}_1(Y)&\to\Z \\
		[a]&\mapsto \#(M^Y(a,\theta)_{\pi,\delta}/\R),
	\end{align*}
	(here $ \delta $ is a small is a small number is the Sobolev space setting) is an  invariant of a homology sphere. Meanwhile, this defines an element $ [\theta^{[s,r]}_Y]\in I^1_{[s,r]}(Y) $ in the filtered instanton cohomology group of $ Y $.
\end{thm}

Indeed, $ [\theta^{[s,r]}_Y]\in I^1_{[s,r]}(Y) $ satisfies the following property

\begin{lemma}[\cite{NST23}, Lemma 3.1]
	For any homology 3-sphere $ Y $ and $ s\in[-\infty,0] $, 
	$$\left\{ r\in(0,\infty]|0=[\theta^{[s,r]}_Y]\in I^{1}_{[s,r]}(Y)\right\}\neq\emptyset.$$
\end{lemma}

Owing to the above lemma, we have

\begin{thm}[\cite{NST23}, Corollary 3.8]
	For $ s\in[-\infty,0] $ and $ Y $ an oriented homology 3-sphere, 
	$$r_s(Y):=\sup\left\{ r\in(0,\infty]\big|0=[\theta^{[s,r]}_Y]\in I^{1}_{[s,r]}(Y)\right\}$$
	is a homology cobordism invariant of $ Y $.
\end{thm}

\subsection{Properties of $ r_s $}
Here we list some properties of $ r_s $ which will be essential for our purpose. 
The first one is an inequality arising from negative definite cobordism.
\begin{prop}[\cite{NST23}, Theorem 1.1(1)]\label{cobordism}
	If there exists a negative definite connected cobordism $ W $ with boundary $ Y_1\amalg -Y_2 $, then the inequality
	$$r_s(Y_2)\leq r_s(Y_1)$$
	holds for any $ s\in[-\infty,0] $.\\
	Moreover, if $ \pi_1(W)=0 $ and $ r_s(Y_1)<\infty $, then the above inequality is strict.
\end{prop}
Not a lot of computations of $ r_s $ are known, but we do know the values for a limited class of Brieskorn homology spheres
\begin{prop}[\cite{NST23}, Corollary 4.6]\label{Brieskorn}
	For any positive coprime $ p,q>1 $ and $ k\in\N $,
	$$r_s(-\Sigma(p,q,pqk-1))=\frac{1}{4pq(pqk-1)},\hsp r_s(\Sigma(p,q,pqk-1))=\infty.$$
	
\end{prop}
The following two properties will be essential when we apply the estimate on $ r_s $ for the surgery description of the branched double cover

\begin{thm}[\cite{NST23}, Corollary 5.6]\label{indep}
	Let $ \{Y_k\}_{k=1}^\infty $ be a sequence of homology 3-spheres satisfying
	$$\infty>r_0(Y_1)>r_0(Y_2)>\cdots,\;\infty=r_0(-Y_1)=r_0(-Y_2)=\cdots.$$
	Then $ \{Y_k\}_{k=1}^\infty $ are linearly independent in $ \Theta^3_\Z $.
\end{thm}

\begin{prop}[\cite{NST23}, Theorem 5.12]\label{qineq}
	For any $ s $ and knot $ K\subset S^3 $, if $ h(S^3_1(K))<0 $, where $h$ is the Fr\o yshov $h$-invariant \cite{Fro02}, then 
	$$\infty> r_s(S^3_1(K))>r_s(S^3_{1/2}(K))>\cdots,\hsp \infty=r_s(-S^3_1(K))=r_s(-S^3_{1/2}(K))=\cdots.$$
\end{prop}

\begin{rmk}
	The assumption of $h(S^3_1(K))<\infty$ is required so that $r_s(S^3_1(K))<\infty$. If by any other way we know that $r_s(S^3_1(K))<\infty$, we can ignore the Fr\o yshov invariant assumption.
\end{rmk}

\subsection{Concordance invariant $\widetilde{s}$ from special cycles}
Fr\o yshov $h$-invariant is notoriously difficult to compute owing to the fact that instanton Floer homology is not well computed except for some Brieskorn spheres \cite{FS90}. $\widetilde{s}$, while difficult to define (we will skip the definition and refer the readers to \cite{DISST22}), it can be computed in many cases owing to the following:
\begin{thm}[\cite{DISST22}, Theorem 1.1, Corollary 1.2]\label{stildeproperty}
	The invariant $\widetilde{s}$ defines a homomorphism $\widetilde{s}:\mc{C}\to\Z$. In particular, if $K$ is alternating, 
	$$\widetilde{s}(K)=-\frac{\sigma(K)}{2}.$$
\end{thm}

The most important application of the $\widetilde{s}$ for our purpose is the following

\begin{thm}[\cite{DISST22}, Theorem 1.5]\label{negcriterion}
	For any knot $K\subset S^3$ with $\widetilde{s}(K)>0$, $h(S^3_1(K))<0$.
\end{thm}

The requirement for $h(S^3_1(K\#J\#K))<0$ in Theorem \ref{main} is certainly more general than the requirement of $\widetilde{s}(K\# J\# K)>0$. However, owing to the fact that $\widetilde{s}$ is a homomorphism and a slice torus invariant \cite{Liv04}, we can easily use it to handle wide classes of patterns, beyond just $\{P_n\}_{n\in\N}$. More examples will be provided in Section 5.

\section{Proof of the main theorems}

\subsection{Proof of the linear independence} With the topological trick in Section 2 and properties of the invariant $r_s$ in Section 3, the proof of the independence of knots in the concordant group is straightforward. 

\textit{Proof of Theorem \ref{main}}: From the algorithm described in Section 2, we know that $$\Sigma_2(Q_n(K))\cong S^3_{1/n}(K\# J\# K).$$ From the assumption that $h(S^3_1(K\# J\# K))<0$, Proposition \ref{qineq}	and Theorem \ref{indep} together implies $\{\Sigma_2(Q_n(K))\}_{n\in\mathbb{N}}$ are linearly independent in $\Theta^3_\Z$. Since $\Sigma_2:\mc{C}\to\Theta^3_{\Z/2}$ is a well-defined homomorphism, this implies $\{Q_n(K)\}_{n\in\mathbb{N}}$ are linearly independent in the concordant group.\qed \\

\textit{Proof of Corollary \ref{maintorus}}: Assumption (2) suffices to imply $h(S^3_1(K\# J\# K))<0$ according to Theorem \ref{negcriterion}.\qed\\

Corollary \ref{whiteheadnble} is immediate from Theorem \ref{maintorus} since $ \widetilde{s}(U)=0 $. Corollary \ref{whitenotslice} also follows because $ r_s(\Sigma_2(P_n(T_{p,q}))<\infty=r_s(S^3) $.



\subsection{Independence under Whitehead $ n $-ble}
To show certain satellite operator is infinite rank, we have to show there exists an infinite family of knots in which its image under the satellite operator is infinite rank in the concordance group.

To proceed in the proof of Theorem \ref{almostrankexpand}, we first recall from for instance \cite{Rol03} that $ \pi_1(S_{1/n}^3(K))=\langle x_1,\cdots,x_l,\lambda|R,\lambda^n x_1=1\rangle $, where $ x_i $ represent meridian of each continuous arc in a diagram of $ K $, $ R $ is the collection of relations coming from the Wirtinger presentation of knot complements and $ \lambda $ is the longitude of the knot $ K $.

\begin{figure*}[h]
	\centering
	\includegraphics[scale=0.9]{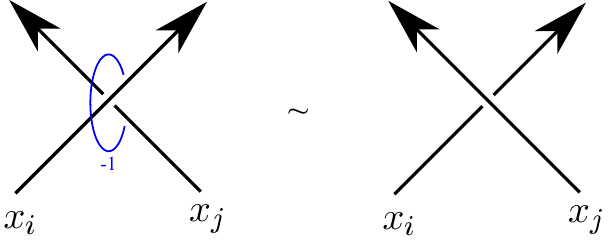}
	\caption{Crossing change by attachment of -1-framed 2-handle}
	\label{crossingchange}
\end{figure*}

If we attach a $-1$-framed 2-handle across of the surgery link as in Figure \ref{crossingchange}, Kirby calculus computation implies this gives a cobordism with two boundary components: one obtained from the surgery over a link $ L $, the other obtained from the surgery over a link $ L' $ with the same framing in which $ L' $ is $ L $ with one crossing changed in sign.

Meanwhile, since $ x_i,x_j $ linked with the core of the 2-handle oppositely, this implies the fundamental group of the cobordism will have a relation of $ x_i\inv{x}_j=1 $.

\begin{figure*}[h]
	\centering
	\includegraphics[scale=0.53]{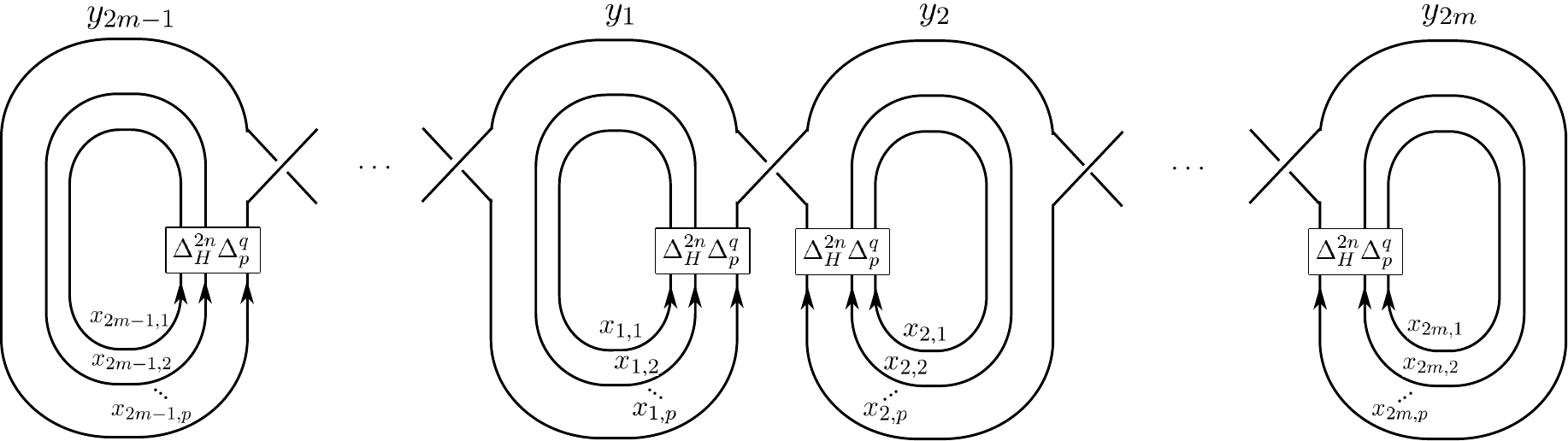}
	\caption{Connected sum of $ 2m $ copies of $ T_{p,q+kp} $}
	\label{2npqtorus}
\end{figure*}

Our goal is to show for $ M=\{(m_k,k)|m_k\leq m_{k+1}\}\subset \N\times\N$, $ \{P_n(mT_{p,q+pk})\}_{(m,k)\in M} $ are linearly independent. Under the Montesinos trick, $ \Sigma_2(P_n(mT_{p,q+kp})) $ is the $ 1/n $-surgery over the knot in Figure \ref{2npqtorus}. Here $ \Delta^q_p $ is the braid word in which its closure is $ T_{p,q} $ while $ \Delta_{H}=\sigma_1(\sigma_2\sigma_1)\cdots(\sigma_{p-1}\cdots\sigma_1) $, where $ \sigma_i $'s are the generators of the $ p $-strand braid group representing a positive crossing represent. $ \Delta_H $ represent a half-twist of $ p $ strands.
In order to prove Theorem \ref{almostrankexpand}, we need the following lemma.

\begin{lemma}\label{generator}
	The group $ \pi_1(S^3\setminus 2mT_{p,q+kp} )$ is normally generated by the meridians of the arcs $ \{x_{i,j}\}_{1\leq i\leq 2m,\;1\leq j\leq p} $ (here we use $ x_{i,j} $ to denote the arc itself and the meridians of arcs interchangingly) in Figure \ref{2npqtorus}.
\end{lemma}

\textit{Proof}: For the arcs lying inside the box of braid words, Lemma 5.21 of \cite{NST23} directly generalize (for the box in the $ i $-th copy of $ T_{p,q+kp} $, all the arcs can be expressed as the conjugation of $ x_{i,j} $'s for $ 1\leq j\leq p $). We only have to show that the arcs $ \{y_i\}_{1\leq i\leq 2m} $ are also just conjugation of $ \{x_{i,j}\} $.

We prove from the right hand most $ y $-arc and show that inductively each of them is generated by $ x_{i,j} $'s. For $ y_{2m} $, it is simply $ x_{2m,p} $. For $ y_{2m-2} $, the Wirtinger relation (see for example \cite{Rol03}) for the positive crossing connecting two copies of $ T_{p,q+kp} $ implies $ y_{2m-2} $ is a conjugation of an arc from the box of braid words by $ y_{2m} $, which is generated by $ x_{i,j} $'s. Inductively, the Wirtinger relation implies all $ y_{i} $ are conjugation by $ \{x_{i,j}\}_{1\leq i\leq 2m,\;1\leq j\leq p} $, proving the claim.\qed \\

With the above lemma in hand, we are ready to construct the necessary cobordism.

\begin{lemma}\label{simpconnected}
	There exists a  cobordism $ W $ which is 
	\begin{enumerate}[(1)]
		\item 
		negative definite;
		\item 
		$ \bd W =S^3_{1/n}(2m_kT_{p,q+kp})\amalg-S^3_{1/n}(2m_{k+1}T_{p,q+(k+1)p})$ for $ m_{k+1}\geq m_k $ and;
		\item 
		$ \pi_1(W)=1 $.
	\end{enumerate}
\end{lemma}

\textit{Proof}: From Lemma \ref{generator}, $ \pi_1(S^3_{1/n}(2m_{k+1}T_{p,q+(k+1)p})) $ is normally generated by $ \{x_{i,j}\}_{1\leq i\leq 2m,\;1\leq j\leq p} $. Note that if we do a positive crossing change for all crossings in one copy of $ \Delta_H $ from each of the $ T_{p,q+(k+1)p} $ by attaching $ -1 $-framed 2-handle as in Figure \ref{crossingchange}, it gives $ \Delta_H^{-1} $ and thus the other end of the cobordism will be $ S^3_{1/n}(2m_{k+1}T_{p,q+kp}) $. This cobordism is negative definite since only $ -1 $-framed 2-handles are attached, proving (1) and (2).  

For (3), note that handle attachments implies $ \pi_1(W) $ is obtained from $ \pi_1(S^3_{1/n}(2m_{k+1}T_{p,q+(k+1)p})) $ by introducing relations  $ x_{i,1}=\cdots=x_{i,p} $ for $ i=1,\cdots,2m $.  Note that $ x_{2m-2,p}=x_{2m,p} $ as they are the same arc in Figure \ref{2npqtorus} and the same argument generalize for $ x_{2m-2,p}=x_{2m-4,p} $ and the subsequent copies of torus knots. This implies $ \pi_1(W) $ is normally generated by \textit{one} element, $ x_{2m,p} $. Without loss of generality, assume the longitude relation from the $ 1/n $-surgery is given by $ \lambda^nx_{2m,p}=1 $. Since $ \lambda $ is in the commutator subgroup of $\pi_1(S^3_{1/n}(2m_{k+1}T_{p,q+(k+1)p}))$ (as $ H_1(S^3_{1/n}(2m_{k+1}T_{p,q+(k+1)p})=0 $), the fact that $ \pi_1(W) $ is normally generated by one element and $ \lambda $ is in the commutator implies $ \lambda=1 $ in $ \pi_1(W) $. So, the longitude relation in $ \pi_1(W) $ implies
$$ x_{2m,p}=(1)^n x_{2m,p}=\lambda^n x_{2m,p}=1,$$
proving (3).

So far, the cobordism we have fixes the number of copies of $ T_{p,q+(k+1)p} $. Note that we can also have a negative definite and simply connected cobordism with boundary being $S^3_{1/n}(2m_kT_{p,q+kp})\amalg-S^3_{1/n}(2m_{k+1}T_{p,q+(k+1)p})$ with $ m_{k+1}>m_k $ just by attaching additional handles on the previously constructed $ W $ which unknot 2 copies of $ T_{p,q+(k+1)p} $.\qed \\

With the above lemma in hand, the proof of \ref{almostrankexpand} is straightforward.\\

\textit{Proof of Theorem \ref{almostrankexpand}}: Lemma \ref{simpconnected} implies the strict inequality of $ r_s $ for the family of homology spheres $ \{\Sigma_2(P_n(m_k T_{p,q+kp}))\}_{(m_k,k)\in M} $ for a set $ M=\{(m_k,k)|m_{k}\leq m_{k+1}\}\subset\N\times\N $ (Proposition \ref{cobordism}) and upperly bounded by $ 1/4pq(pqn-1) $ (there exists a negative definite cobordism arising from crossing change to $ S^{3}_{1/n}(T_{p,q})\cong-\Sigma(p,q,pqn-1) $ and apply Proposition \ref{Brieskorn}).\\
Meanwhile, as $ \Sigma_2(P_n(m_kT_{p,q+kp})) $ can be represented as a $ 1/n $-surgery over a knot, Proposition \ref{qineq} implies $ r_s(-\Sigma_2(P_n(m_kT_{p,q+kp})))=\infty $. Putting everying together, this implies the family $ \{\Sigma_2(P_n(m_k T_{p,q+kp}))\}_{(m_k,k)\in M} $ is linearly independent in $ \Theta^3_\Z $, thus concluding the linear independence of $ \{P_n(m_k T_{p,q+kp})\}_{(m_k,k)\in M}\subset\mathcal{C} $.\qed

\section{Examples for Theorem \ref{main} with nontrivial $J$ and non-torus companion}

In this last section, we devote to provide some examples in which Theorem \ref{main} is applicable but the Montesinos friend $J$ is nontrivial and the companion knot $K$ is not merely a torus knot.

\begin{figure*}[h]
	\centering
	\includegraphics[scale=0.8]{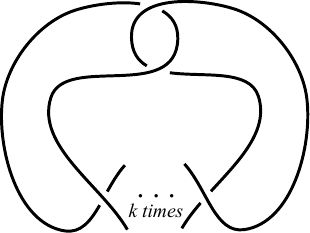}
	\caption{Diagram of the twist knot $K_n$. Here $ k $ denotes the number of half-twists in the horizontal direction. If $ k $ is negative, it means $ |k| $ negative half-twists.}
	\label{twistknot}
\end{figure*}

\begin{prop}
	Let $\{Q_n\}_{n\in\N}$ be the family of pattern in Figure \ref{realization} with tangles in the green circle which corresponds to $3+1/n$ in the standard enumeration, then for a fixed twist knot $K_m$ with $m\in\Z$ (see Figure \ref{twistknot}) $\{Q_n(K_m)\}_{n\in\N}$ are linearly independent.
\end{prop}
\textit{Proof}: From the algorithm described in Section 2, $\Sigma_2(Q_n(U))\cong S^3_{1/n}(T_{2,3})$. From Theorem \ref{stildeproperty}, this implies $\widetilde{s}(K_m\#T_{2,3}\# K_m)=1-\sigma(K_m)>0$ since it can be computed that 
$$\sigma(K_m)=\begin{cases}
	-2 & \text{if }m\text{ is odd}\\
	0 & \text{if }m\text{ is even}
\end{cases}.$$
Hence, from Corollary \ref{maintorus}, this implies $\{Q_n(K_m)\}_{n\in\N}$ are linearly independent.\qed

	\printbibliography

\end{document}